\documentstyle[12pt,twoside,epsfig]{article} 
\begin{document} 
\null 
\vskip 1 cm 
\centerline {\bf COMPUTER PROGRAMS FOR KNOT TABULATION}
\bigskip \centerline
{Charilaos Aneziris} \par 
\bigskip
\centerline {\it Institut f\"ur Hochenergienphysik, Zeuthen}
\centerline {\it D.E.SY. (Deutsches Elektronen-SYnchrotron)}
\centerline {\it Platanenallee 6, Zeuthen 15738} 
\centerline {\it GERMANY} \par
\bigskip \centerline {\bf Abstract} \par \bigskip
     While the problem of knot classification is far from solved, it is 
     possible to create computer programs that can be used to tabulate knots 
     up to a desired degree of complexity. Here we discuss the main ideas on 
     which such programs can be based. We also present the actual results 
     obtained after running a computer program on which knots are denoted 
     through regular projections. \par \bigskip 
\centerline {\bf 1. \quad INTRODUCTION} \par \bigskip  
     In order to develop an algorithm to tabulate knots, one first needs to 
     introduce an appropriate notation, so that knots may be presented through 
     sequences of numbers. Such a notation may only be valid if each acceptable 
     presentation uniquely determines a knot up to ambient isotopy. One then 
     continues by obtaining the conditions that have to be satisfied by some 
     sequence of numbers, in order for the sequence to be acceptable and thus 
     denote an actual knot. This problem is far from trivial, as shall be made 
     clear later in the paper. 
     \par \medskip 
     The next step is to identify the sequences that denote ambient isotopic 
     knots, through the introduction of appropriate equivalence moves. If two 
     or more number sequences are related through such moves, all sequences 
     except the one appearing first are eliminated and do not appear at the 
     output, and thus the corresponding knot type is tabulated exactly once.
     \par \medskip 
     There is no upper bound however to the number of equivalence moves that 
     may be needed to connect two sequences denoting ambient isotopic knots. 
     Therefore if two sequences are not connected after some finite number of 
     moves, this does not imply that the corresponding knots are inequivalent. 
     In order to show inequivalence one needs to establish appropriate 
     invariants. Each such invariant is a function on the set of notations, 
     defined so that any sequences connected through equivalence moves are 
     mapped to identical values. Therefore if one or more invariants map two 
     sequences to distinct values, the corresponding knots are inequivalent.
     \par \medskip 
      Ideally one may run the program long enough, so that each knot up to some 
      desired complexity has either been shown ambient isotopic to a preceding 
      knot, or shown inequivalent to all preceding knots. In the former case 
      the knot is eliminated, while in the latter it is tabulated. Due however 
      to memory and time limitations, there are certain limits up to which one 
      may carry such a process. \par \medskip 
     In Section 2 we discuss two methods of denoting knots. One is based on the 
 {\bf braid word}, and the other on presenting knots as closed Self-Avoiding 
     Walks (SAW's) on a cubic lattice. In Section 3 we present and discuss in 
     detail the {\bf Dowker notation}, which is based on studying knots through 
     their regular projections. Finally in Section 4 we present the results we 
     obtained after running a computer program based on the Dowker notation, 
     and compare them to results obtained by others (Ref. 1). \par \bigskip 
\centerline {\bf 2. \quad EXAMPLES OF KNOT NOTATIONS} \par \bigskip
\centerline {$\underline {\bf The \quad Braid \quad Word}$} \par \bigskip
  According to the {\bf Alexander Theorem}, every knot and link may be obtained 
   through the closure of appropriate braids (Ref. 2,3). Let a knot $K$ be the 
     closure of a braid where  
     $b= \prod _{i=1} ^k \sigma _i ^{\epsilon _i}$ and $\sigma _i \in   
     \{ 1,2,...,n-1 \}$.  
     is a generator of the braid group $B_n$. One may 
     thus denote $K$ by the sequence $a_0,a_1,a_2,...,a_k$ where
     $a_0=n$ and $a_i=\sigma _i \epsilon _i$. 
     Once this sequence is given, one 
     may uniquely determine the braid $b$ and its closure $K$. 
     \medskip In general such 
     sequences yield multicomponent links, and one thus needs a criterion to 
     distinguish the sequences that denote knots from the ones that denote 
 links. In order to do so, one replaces the braid generators $\sigma _i$ with 
 permutations $p_i=(i,i+1) \in S_k$. The permutation  
 $P=\prod _{i=1} ^k p_i$ belongs to a conjugacy class of $S_k$ characterized by
 a partition $\rho _1 \geq \rho _2 \geq ... \geq \rho _k > 0$  
     of $k$. The braid closure of ${\prod _{i=1} ^k} \sigma _i ^{\epsilon _i}$
     is a $q$-component link. If $q=1$, such a link is a knot. 
     \par \medskip Two distinct sequences 
    denote the same knot if and only if they either correspond to identical 
    braids, or to braids connected through equivalence moves. The first 
    possibility occurs by inserting or removing trivial factors 
    $\sigma _l ^{\epsilon } \sigma _l ^{- \epsilon }$, by replacing
    $(\sigma _i \sigma _j ) ^{\epsilon }$  
    with $(\sigma _j \sigma _i ) ^{\epsilon }$  
    for $|i-j|>1$, or by replacing  
	$(\sigma _i \sigma _{i+1} \sigma _i) ^{\epsilon }$  
with $(\sigma _{i+1} \sigma _i \sigma _{i+1}) ^{\epsilon }$.   
    \par \medskip Two distinct braids yield 
    ambient isotopic knots if they are connected by a series of the following 
    moves (Ref. 4). \par \medskip 
    $\alpha ) \qquad b \leftrightarrow {b' b b'} ^{-1}$
where $b,b' \in B_k$ \par \smallskip $\beta ) \qquad b \leftrightarrow 
b \sigma _k ^{\epsilon }$
where $b \in B_k, \quad \epsilon = \pm 1$  
and $b \sigma _k ^{\epsilon } \in B_{k+1}$ \par \medskip 
Recently (Ref. 3) it has been shown that these two kinds of 
moves, known as the {\bf Markov} moves, can be generalized and replaced by the 
$L$-moves \par \smallskip
$$b=b_1 b_2 \leftrightarrow \sigma b_1 ' \sigma ' \sigma _k ^{\pm 1} (\sigma ' )
^{-1} b_2 ' \sigma ^{-1}$$ \par \smallskip \noindent   
where $\sigma = \sigma _i ^{\epsilon } ... \sigma _k ^{\epsilon }$,  
$\sigma ' = \sigma _{i-1} ^{-\epsilon } ... \sigma _{k-1} ^{-\epsilon }$
and $b_1 '$, $b_2 '$ 
are derived from $b_1$, $b_2$ by replacing $\sigma _j$ with $\sigma _{j+1}$ 
$\forall j \geq i$. \par \medskip Since 
there is no upper bound in the number of the Markov or L-moves that may be 
needed, one has to establish appropriate invariants in order to distinguish 
braids yielding inequivalent knots. Such invariants may be obtained through the 
{\bf skein} relations \par \medskip {\begin {figure}[h]  
$$A f(\epsfysize=0.5 cm \epsffile {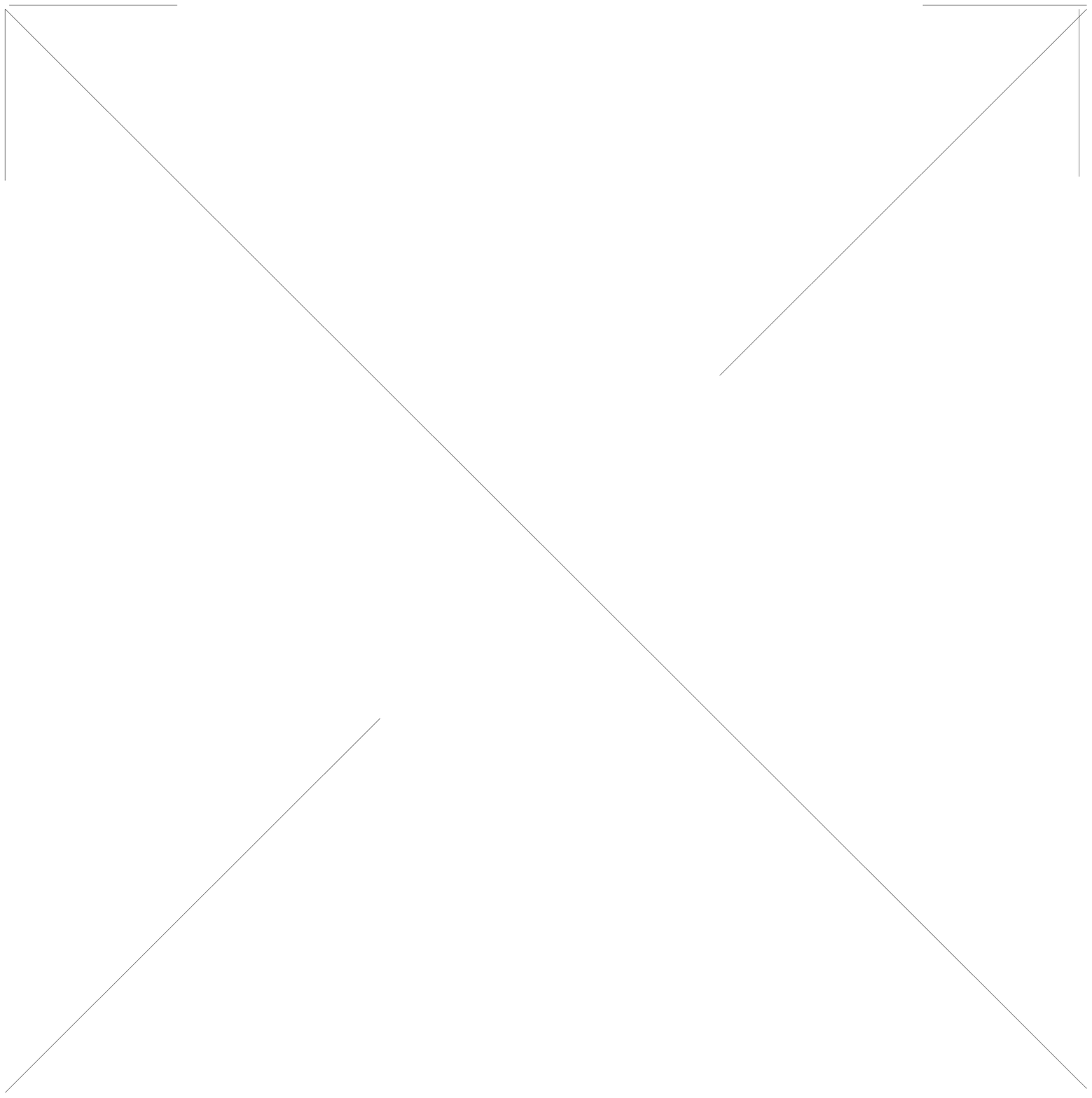})  
+ B f(\epsfysize=0.5 cm \epsffile {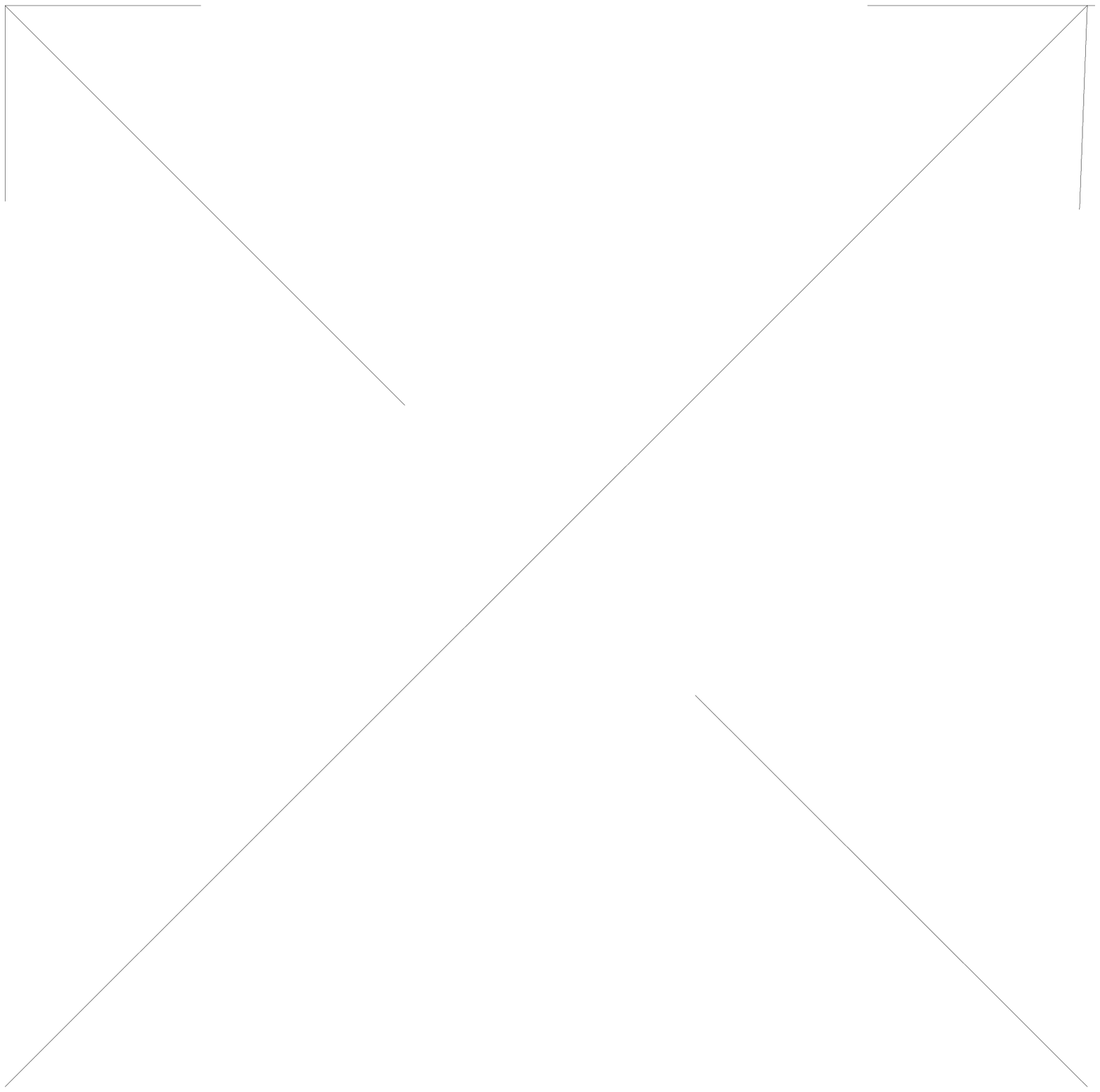})  
= C f(\epsfysize=0.5 cm \epsffile {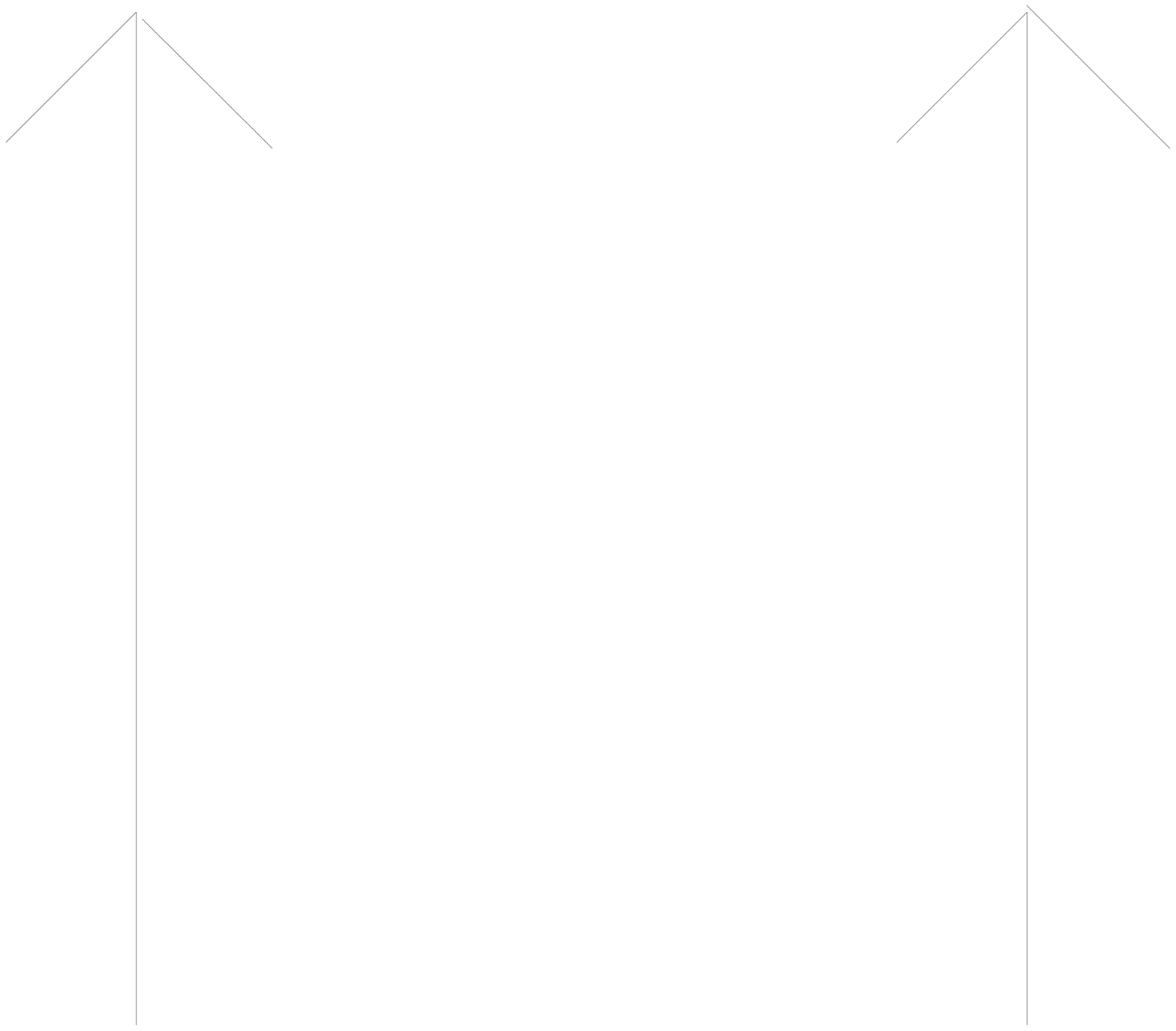}) $$  
\end {figure}}   
	  \par \medskip \noindent  
	  from which one deduces that 
	  $Af(b_1 \sigma _i b_2 ) + Bf(b_1 \sigma _i ^{-1} b_2)=cf(b_1 b_2 )$. 
	  Examples of 
	  such invariants are the 
	  {\bf Alexander-Conway}, the {\bf Jones} and the {\bf homflypt} 
	  polynomials (Ref. 5). \par  
	  \medskip Such invariants however do not always 
	  distinguish inequivalent knots. In addition, knots expressed as braid 
	  closures may not appear in a minimum crossing number projection, as 
	  the knot 5-2 in Rolfsen's tables demonstrates (Ref. 6). Therefore the 
	  braid word notation becomes inefficient when the crossing number is 
	  large. \par \bigskip
\centerline {$\underline {\bf Closed \quad SAW's \quad on \quad a \quad Cubic 
\quad Lattice}$} \par \bigskip
     Let $\{ (m,n,l)/m,n,l \in Z\}$  
     define a cubic lattice in $R^3$. By connecting neighboring points on the 
lattice, one obtains a {\bf walk}. If the initial and final points coincide, 
   and no other point appears more than once, the walk is called closed 
Self-Avoiding, and defines a knot. A {\bf tame} knot is defined as a knot that 
     is ambient isotopic to a closed Self-Avoiding Walk. Therefore the problem 
of tabulating (tame) knots can be reduced to tabulating closed SAW's. 
\par \medskip
      Let $\vec {r_0}, \vec {r_1}, \vec {r_2}, ... , \vec {r_n} = \vec {r_0}$  
      the points on the cubic lattice defining the walk, such that 
      $|\vec {r_i} - \vec {r_{i-1}}| = 1 \forall i \in \{ 1,2,...,n \}$. The 
      unit vectors $\vec {v_i} = \vec {r_i} - \vec {r_{i-1}}$  
      take values in $\{ \pm \vec {i}, \pm \vec {j}, \pm \vec {k} \}$ 
      which may be identified with the numbers 
      $\{ 1,2,3,-1,-2,-3 \}$. One may thus denote a closed SAW and its ambient 
      isotopy class through a sequence of $n$ numbers taking values in 
 $\{ 1,2,3,-1,-2,-3 \}$ and indicating the vectors $\vec {v_i}$. 
 In order for some random sequence $a_1,a_2,...,a_n$, where $a_i \in
\{ \pm 1, \pm 2, \pm 3 \}$, to denote a closed SAW, two conditions have to be 
      satisfied. First, the sum of the corresponding unit vectors  
      $\sum _{i=1} ^n \vec {v_i} = \vec 0$ in order for 
 the walk to be closed. Second, there should exist no two numbers $i$, $j$, 
where $i<j$ and $\{ i,j \} \neq \{ 0,n \}$, such that
$\sum _{k=i+1} ^j \vec {v_k} = \vec 0$, in order for the walk to be 
Self-Avoiding. \par \medskip 
Distinct notations $a_i$, $b_i$ yield identical walks if there is a number  
$c \in \{ 1,2,...,n-1 \}$ such that
$\forall i \in \{ 1,2,...,n-c \}$, $b_{i+c}=a_i$ and 
$\forall i \in \{ n-c+1,...,n-1,n \}$, $b_{i+c-n}=a_i$,
since a switch from  to  merely relocates the end points of 
     the walk, and is thus a reparametrization of the walk. 
     \par \medskip Distinct walks 
     yield ambient isotopic knots if they are related through the equivalence 
     moves \par \medskip {\begin {figure}[h] \hspace {1.2 cm} {\bf I)} \qquad  
\epsfysize=0.5 cm \epsffile {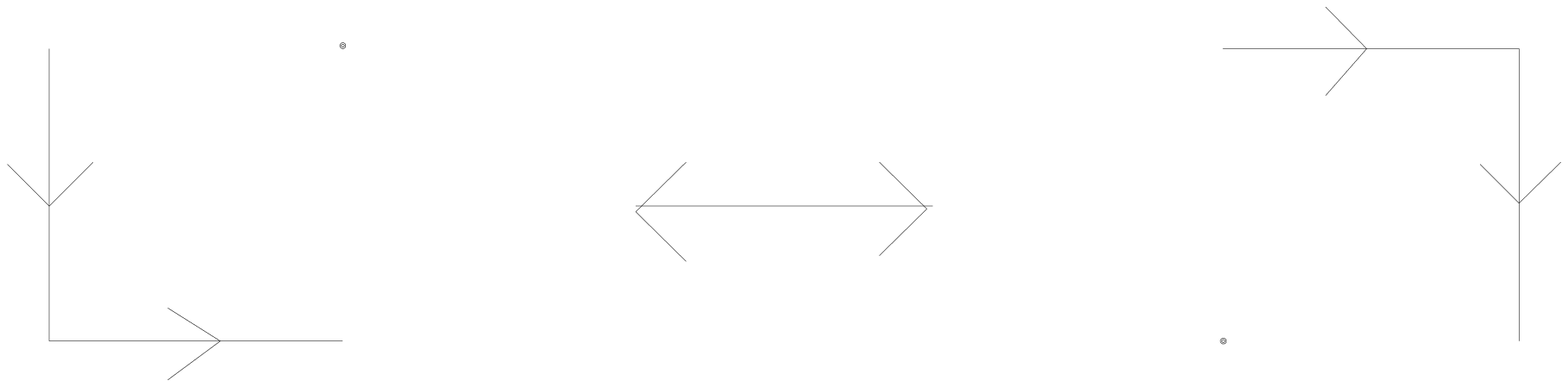} \qquad  
$b_i=a_{i+1}, \qquad b_{i+1}=a_i, \qquad j \neq i,i+1 \Rightarrow
     b_j=a_j$ \end {figure}} {\begin {figure}[h] \hspace {1 cm} {\bf II)}   
    \qquad \epsfysize=0.5 cm \epsffile {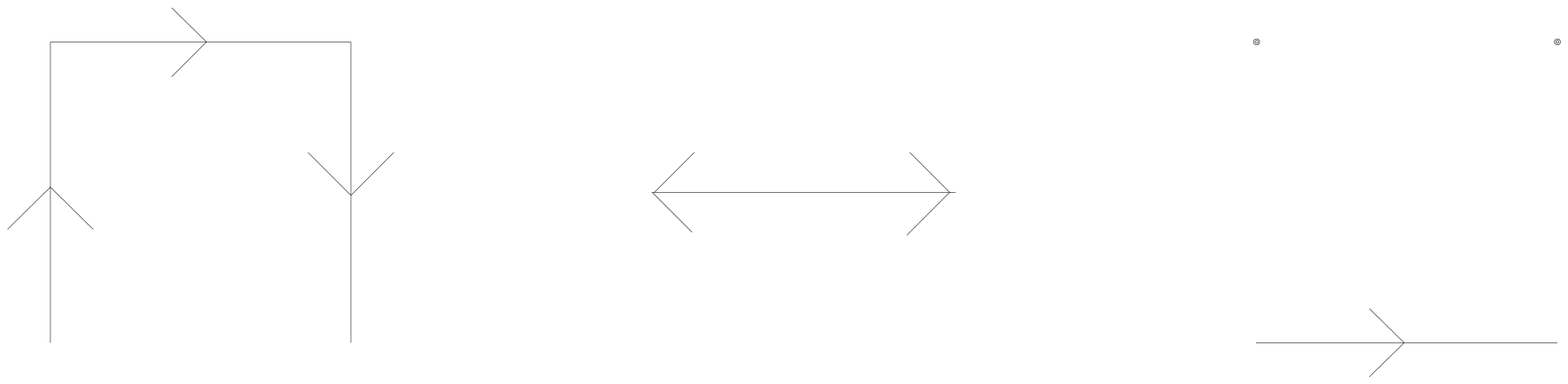} \qquad   
$j \leq i-1 \Rightarrow b_j=a_j, \qquad b_i=a_{i+1}, \qquad
     j \geq i \Rightarrow b_j=a_{j+2}$ \end {figure}}   
     \hskip 4.8 cm (Possible only if $a_{i+2}=-a_i$) 
     \par \bigskip While under such moves a closed SAW 
     remains closed, it may not remain Self-Avoiding. One thus has to ensure 
    that any new lattice points to be occupied after a move, were unoccupied 
    before, in order for the move to be allowed (Ref. 7). 
    \par \medskip As was the case with 
    the braid word notation, there is no upper bound in the number of moves 
    needed to connect two ambient isotopic walks, and one thus needs to 
    establish invariants. In order to do so, one first performs the necessary 
    equivalence moves so that all segments along one of the six directions 
    belong to the surface of a rectangular parallelepiped that encloses the 
    walk. One then removes them and converts the closed walk to a collection of 
    open walks, which are considered as the strings of a braid whose closure is 
    the knot. It is possible from these open walks to obtain the braid word, 
    and continue as before by calculating the skein invariants (Ref. 8).
    \par \medskip 
      In addition however to the problem of inequivalent knots with identical 
      skein invariants, running a program based on the cubic lattice notation 
      is extremely timing consuming. The simplest non-trivial knot, for 
      example, the trefoil, needs a minimum length of 24 in order to be 
      represented as a closed SAW on a cubic lattice; therefore the program 
      will spend an enormous amount of time merely showing that all walks of 
      length shorter than 24 are trivial (Ref. 9). In contrast, one may obtain 
      it very fast through other notations, since its braid word is  
      $\sigma _1 \in B_2$ and its 
      crossing number is 3. \par \bigskip 
      \centerline {\bf 3. \quad THE DOWKER NOTATION} \par \bigskip 
\centerline {$\underline {\bf Denoting \quad Regular \quad Projections}$}
      \par \bigskip 
Let a knot $K$ and $P(K)$ a two-dimensional projection of $K$. $P(K)$ is called 
      a {\bf regular} projection of $K$ if none of its multiple points has a 
      neighborhood that looks as follows. \par \smallskip 
{\begin {figure}[h] \centerline {\bf {I)} \hspace {0.5 cm} \epsfysize = 0.7 cm  
\epsffile {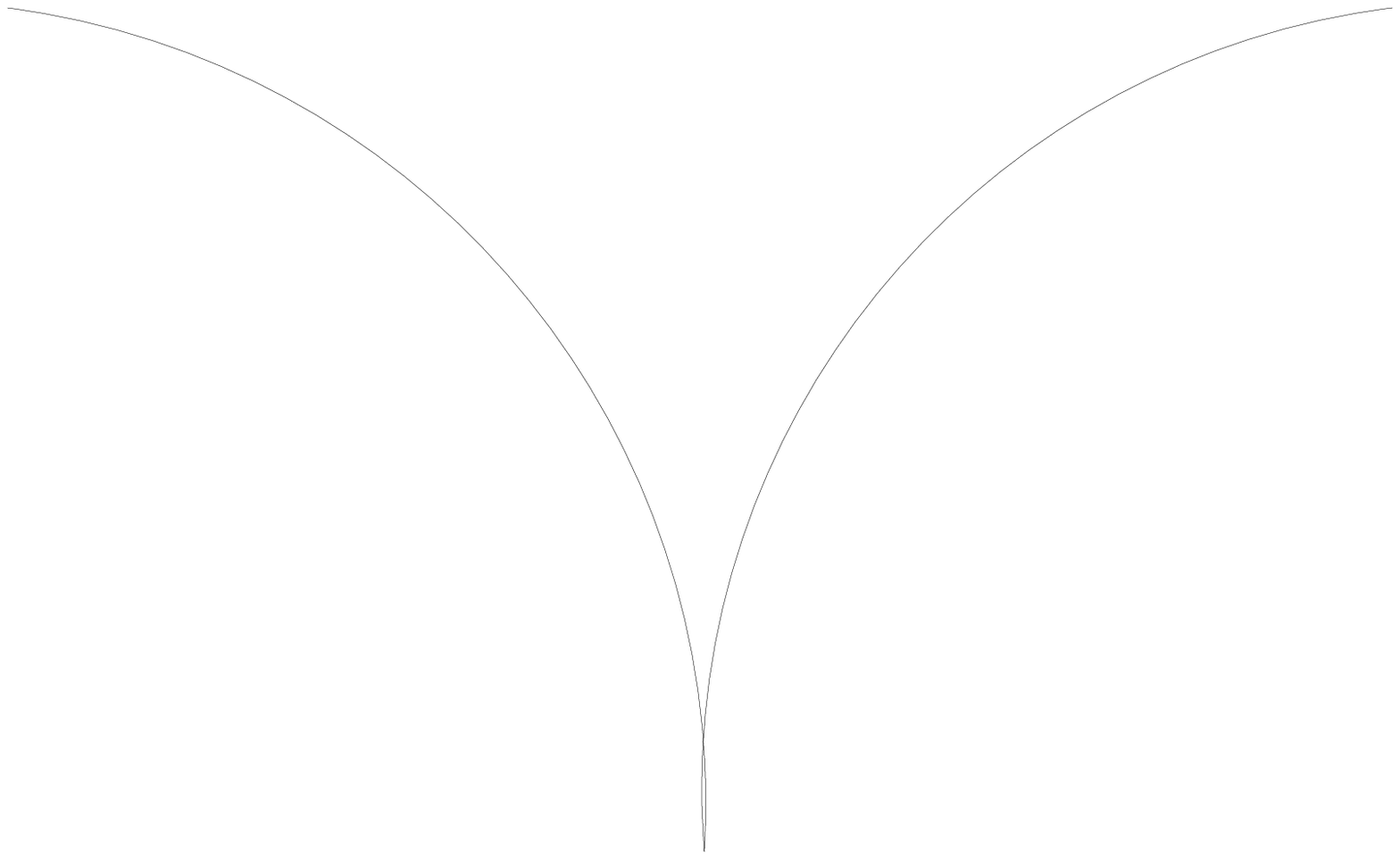} \hspace {0.5 cm} \bf {II)} 
\hspace {0.5 cm} \epsfysize = 0.7 cm  
      \epsffile {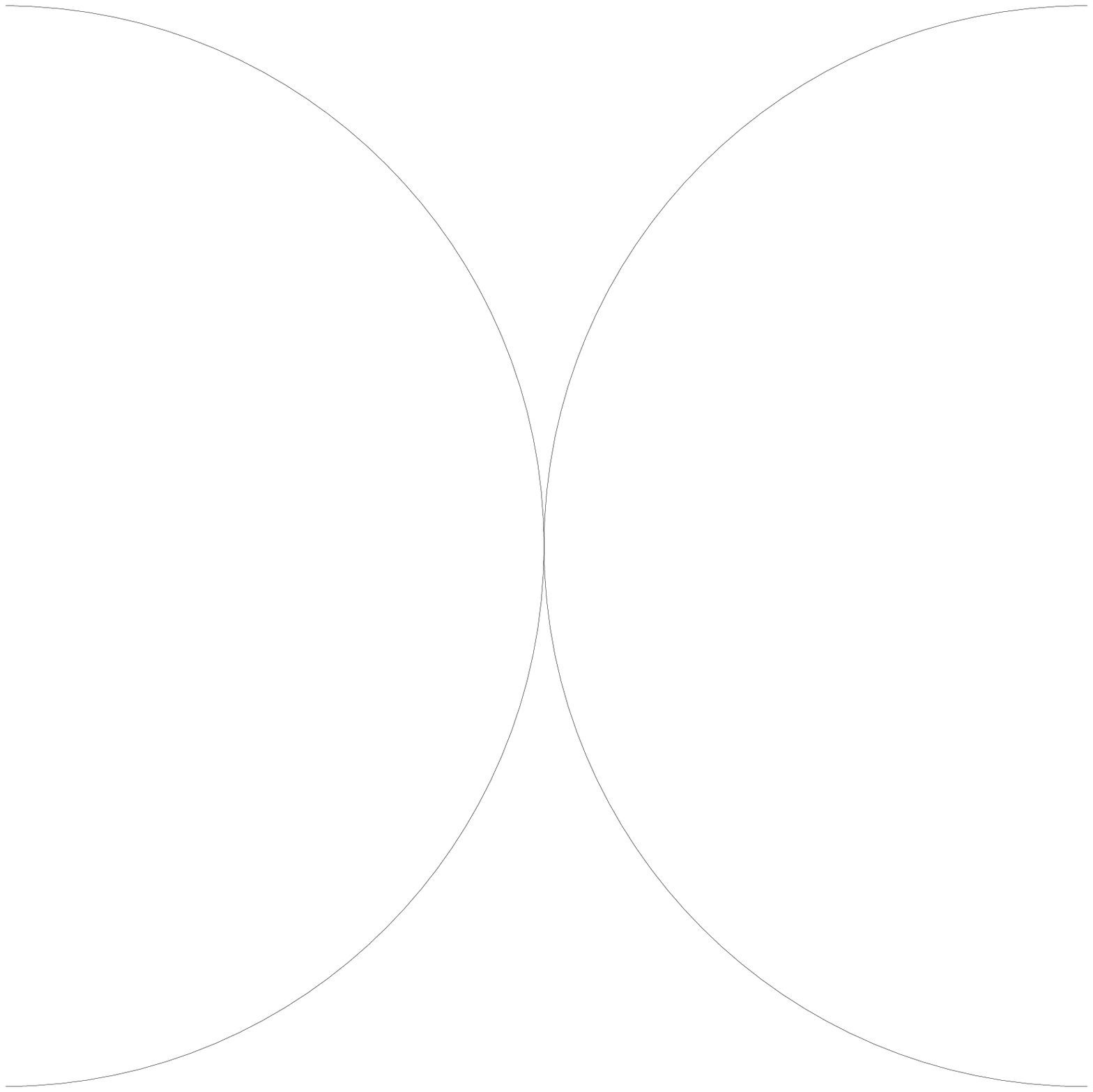} \hspace {0.5 cm} \bf {III)} 
\hspace {0.5 cm} \epsfysize = 0.7 cm \epsffile {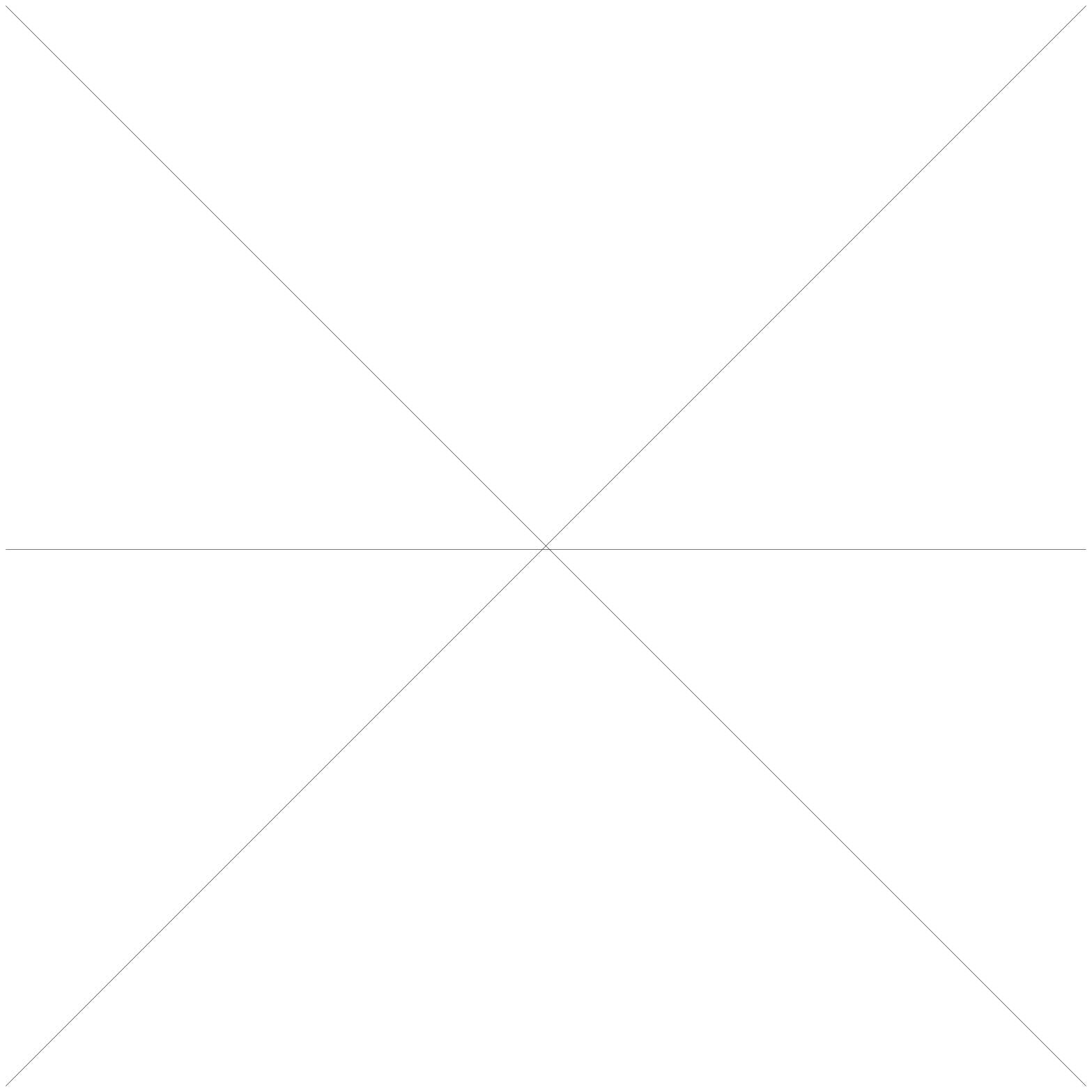}} 
\hspace {0.5 cm}     
\end {figure}} \par \smallskip          
      In other words, let  
      $\vec f (s)$ be the continuous one-to-one function from $S^1$ to $S^3$  
 that defines the knot, and let $\vec t$ be a vector normal to the projection 
      surface. In order for the projection to be normal, the following 
      conditions must be met. \par \medskip     
      I) $\forall s \in S^1 \exists \vec f'(s) = {d \vec f(s) \over ds} :
      \vec f'(s) \times \vec t \neq \vec 0$ \par \smallskip 
      II) $\forall s_1 \neq s_2 \in S^1: \Big( \vec f(s_2) - \vec f(s_1) \Big)
      \times \vec t = \vec 0, \quad \Big( \vec f'(s_2) - \vec f'(s_1) \Big)
      \times \vec t \neq \vec 0$ \par \smallskip
      III) $s_1, s_2, s_3 \in S^1, s_1 \neq s_2 \neq s_3 \neq s_1 \quad
      \Rightarrow \quad |\Big( \vec f(s_2) - \vec f(s_1) \Big) \times \vec t |
   + |\Big( \vec f(s_3) - \vec f(s_1) \Big) \times \vec t| > \vec 0$ 
   \par \medskip
      In addition, in regular projections one distinguishes overcrossings from 
      undercrossings.    
      \par \medskip Regular projections are denoted as follows. Starting 
      from some point and moving along the knot, one assigns successive 
  natural numbers $1$, $2$, ..., $2n$ to the crossing points according to the 
      order these points are met. Each crossing is thus assigned two such 
numbers, $a_0$ for the overcrossing point and $a_0$ for the undercrossing. The
set of the ordered pairs $(a_0,a_u)$ denotes the projection. \par     
\medskip If one considers 
      mirror symmetric knots equivalent and disregards connected sums, given 
      an appropriate set of numbers one can uniquely determine the projection.
\par \bigskip \centerline {$\underline {\bf Drawable \quad - \quad Undrawable 
\quad Sets}$} \bigskip 
     Let n be the crossing number. Each number in 
     $\{ 1,2,...,2n \}$ appears exactly 
     once in one of the $n$ pairs. Not all such sets however denote actual 
     projections, a simple counterexample being $\{ (1,3),(2,4) \}$. 
     \par \medskip
{\begin {figure}[h] \centerline {\epsfysize = 4 cm \epsffile {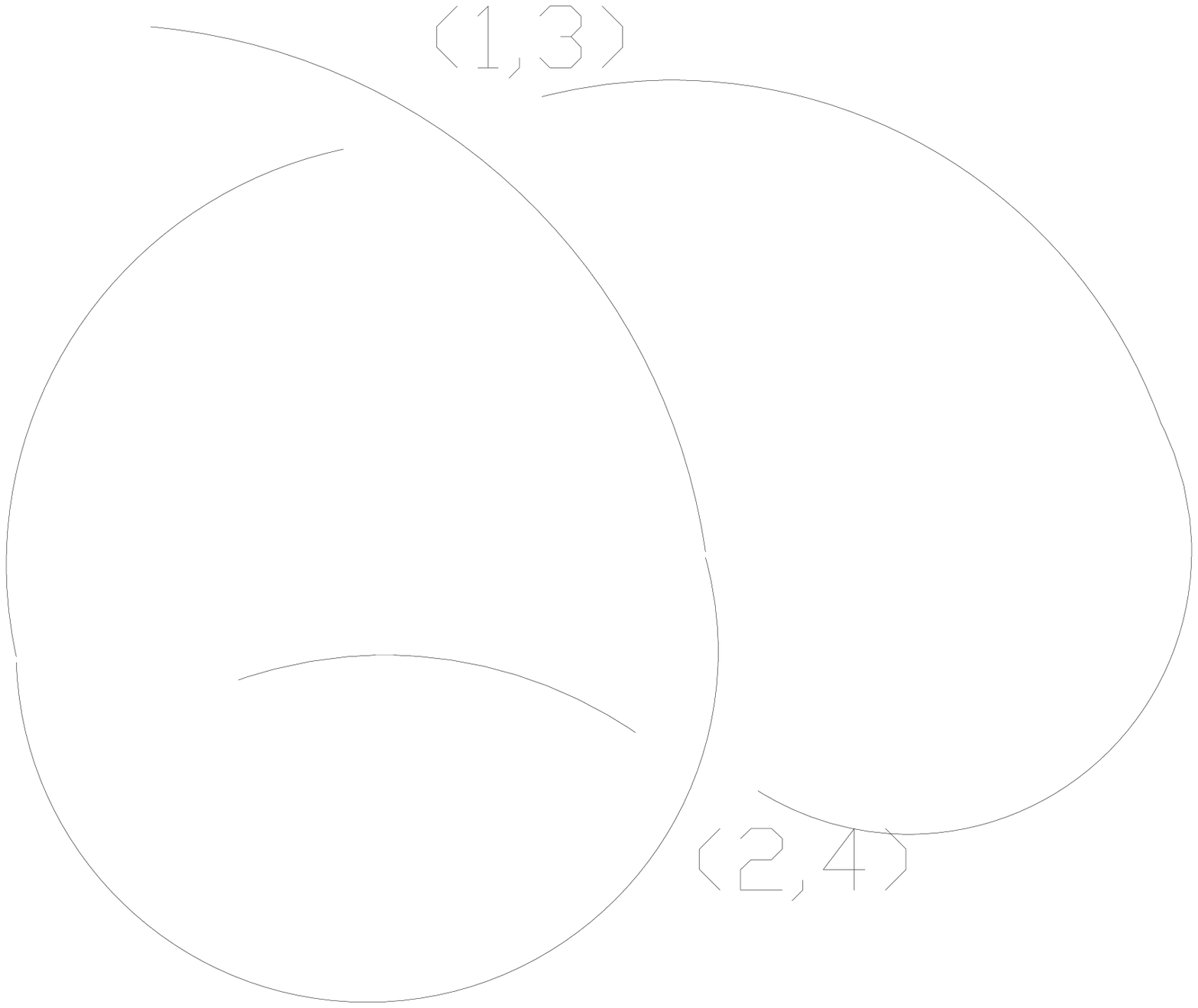}}     
    \end {figure}}  
\par \medskip
     This set, as one notices from the figure above, yields no actual 
     projection and may thus be considered ``undrawable". The reason is the 
     following. Since $(1,3)$ belongs to the set, the segments 
     $1-2-3$ and $3-4-1$ 
     are loops. If $\{ (1,3),(2,4) \}$ were drawable, the loop 
     $1-2-3$ would have 
     divided the plane into two disjoint segments according to the {\bf Jordan 
     Curve Theorem} (Ref. 10). Therefore the loop 
     $3-4-1$ would either be tangent 
     to $1-2-3$, or would have to intersect it at an even number of points, in 
     order to ``enter" as many times as it ``exits" (vertices do not count). 
These two loops are nowhere tangent, since they share no common segment, and  
they intersect at just one point, namely at $\{ (2,4) \}$. Therefore the set 
   is undrawable. 
\par \medskip 
      The general condition in order for a set to be drawable is the 
      following. Any two loops defined by such a set must either share one or 
      more segments, or intersect at an even number of points, vertices not 
      counting.
\par \medskip 
     A corollary of this condition is that odd numbers must always be paired 
     to even. If that were not the case, there would be some pair $(i,j)$ 
     belonging to the set such that $|j-i|$ would be even. Therefore the loops 
$i,i+1,...,j-1,j$ and $j,j+1,...,2n-1,2n,1,...,i-1,i$ which share no common 
     segment, would intersect at an even number of points and thus violate the 
     condition. This corollary, however, while being necessary, is not 
sufficient, as the counterexample $\{ (1,4),(3,6),(5,8),(7,10),(9,2) \}$ shown 
below  demonstrates. In this case, one notice that the loops $1-2-3-4$ and 
$5-6-7-8$ violate the drawability condition. \par \medskip
{\begin {figure}[h] \centerline {\epsfysize = 4 cm \epsffile {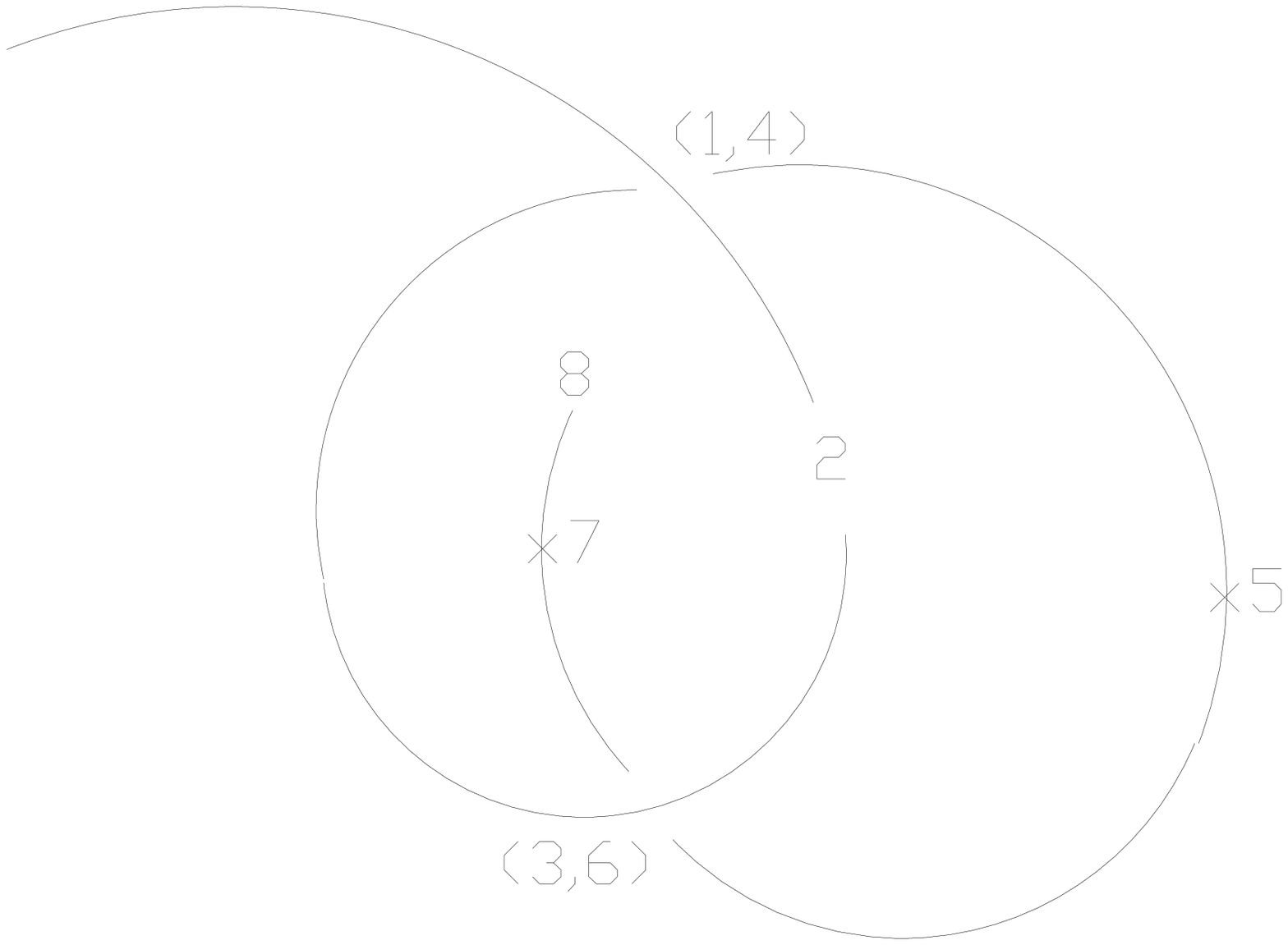}}   
    \end {figure}}  
\par \medskip
     In order to check whether one set is drawable, one has to consider all 
     possible loops, which are $3^n$, since each crossing may or may not be a 
     vertex of the loop. In addition, if it is a vertex, the loop may turn 
     clockwise or counterclockwise. The number of pairs of loops is thus 
     $3^{2n}$. If 
     one was to check all such pairs of loops for each possible set, the CPU 
     time needed would be enormous and grow exponentially. This is not 
necessary, however, since the sets may be divided into classes of {\bf equal 
     drawability}. Sets belonging to one class are either all drawable, or all 
     undrawable, and thus one has to check the pairs of loops in merely one 
     set of each class. For $n=14$, for example, there are $14!$  (about one 
     hundred billion) possible sets. (While checking drawability, one need not 
     consider the ordering in each pair of numbers. When this is also taken 
     into account, the number of possible sets becomes  
     $2^{14}=16,384$ times larger). These 
     however can be divided into less than a thousand classes of equal 
     drawability, and thus the time needed is significantly reduced.     
     \par \medskip Two 
     projections are ``equally drawable" and belong to the same class, if they 
     can be connected through the following moves. \par \medskip 
{\begin {figure}[h] \hspace {2 cm} {\bf I)} \hspace {3 cm} \epsfysize = 0.9 cm 
\epsffile {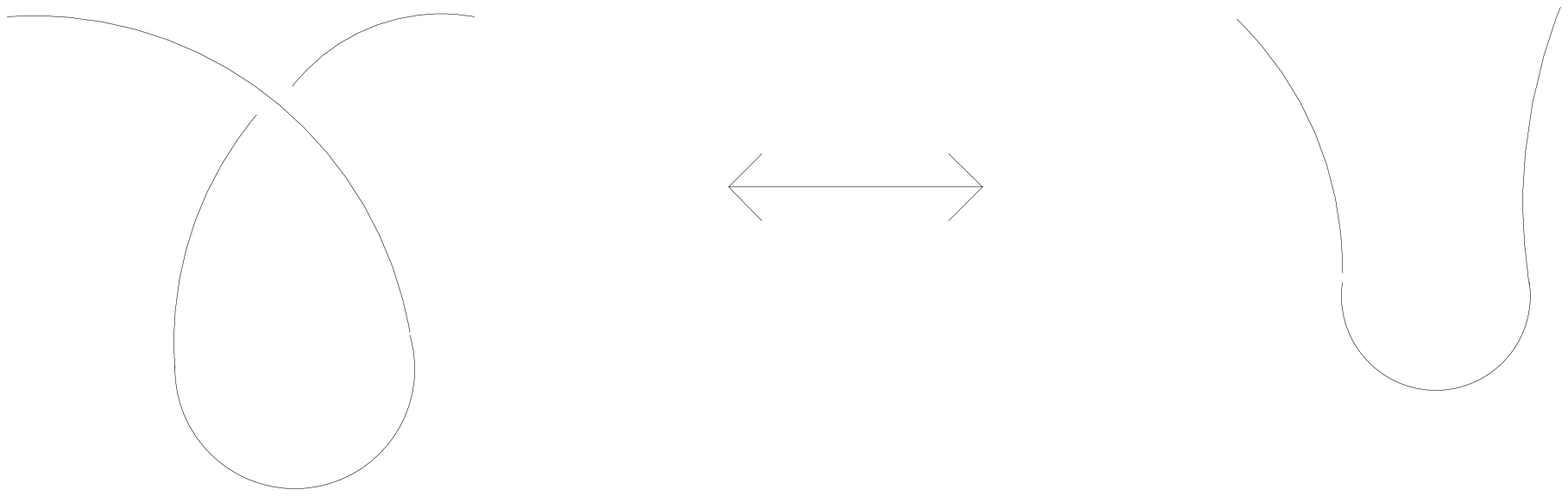} \end {figure}}  
\par \smallskip Add or remove a pair $(i,i+1)$, and increase or 
     decrease numbers larger or equal to $i$ by $2$.    
     \par \smallskip
{\begin {figure}[h] \hspace {1.8 cm} {\bf II)} \hspace {3 cm} 
\epsfysize = 0.9 cm 
\epsffile {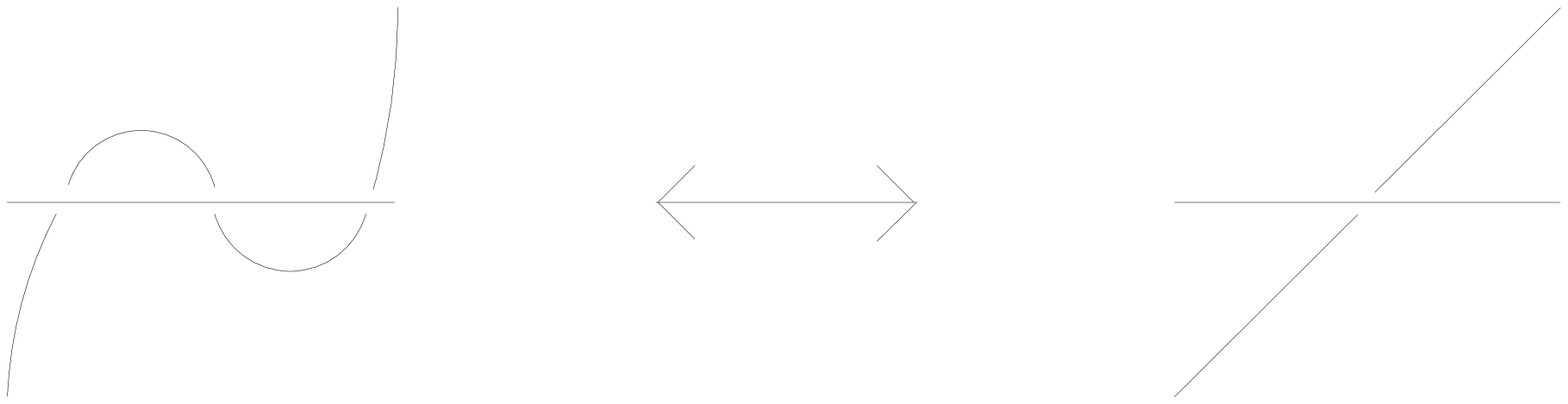} \end {figure}}  
\par \smallskip If a set 
     contains pairs $(i,j),(i\pm 1,j\pm 1),(i\pm 2,j\pm 2)$ 
     remove two consecutive pairs, decrease numbers between $i$ and $j$ by $2$, 
numbers larger than $j$ by $4$. \par \smallskip   
{\begin {figure}[h] \hspace {1.6 cm} {\bf III)} 
\hspace {3 cm} \epsfysize = 0.9 cm 
\epsffile {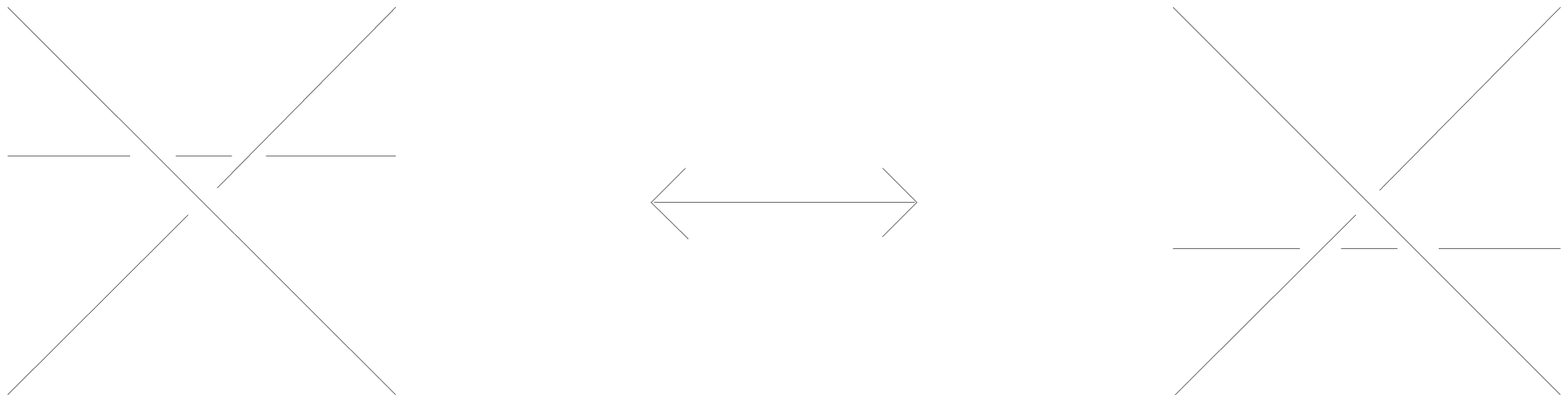} \end {figure}}  
\par \smallskip Substitute pairs $(i,j),(i',k),(j',k')$ with 
  pairs $(i,k'),(i',j'),(j,k)$ where $|i'-i|=|j'-j|=|k'-k|=1$. \par \bigskip  
     \centerline {$\underline {\bf Equivalence \quad Moves}$} \par \bigskip 
     According to the {\bf Reidemeister Theorem} 
     (Ref. 11), two regular knot projections belong to ambient isotopic knots 
     if and only if they are connected through a series of the following moves.
    \eject  
{\begin {figure}[h] \hspace {2 cm} {\bf I)} \hspace {3 cm} \epsfysize = 0.9 cm 
\epsffile {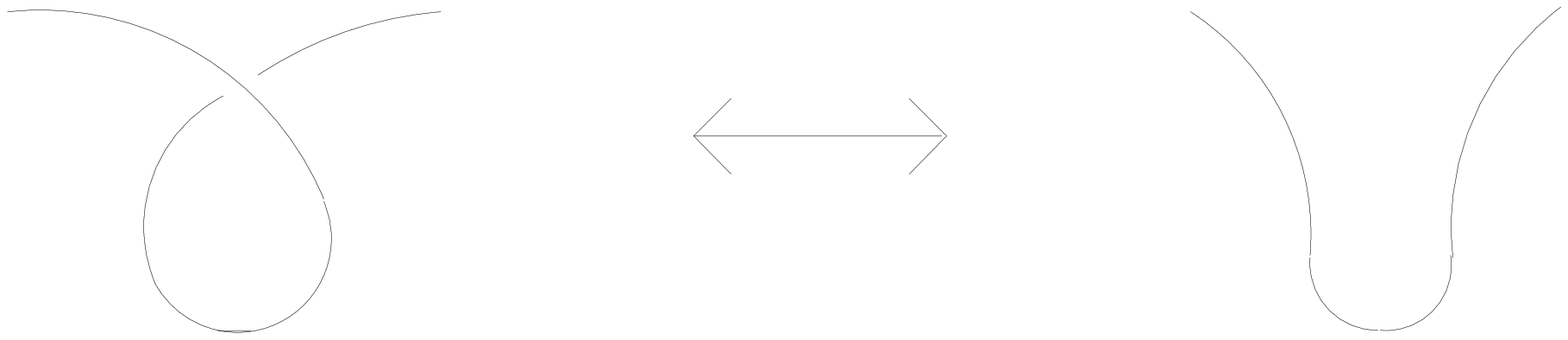} \end {figure}} \par  
{\begin {figure}[h] \hspace {1.8 cm} {\bf II)} \hspace {3 cm} 
\epsfysize = 0.9 cm 
\epsffile {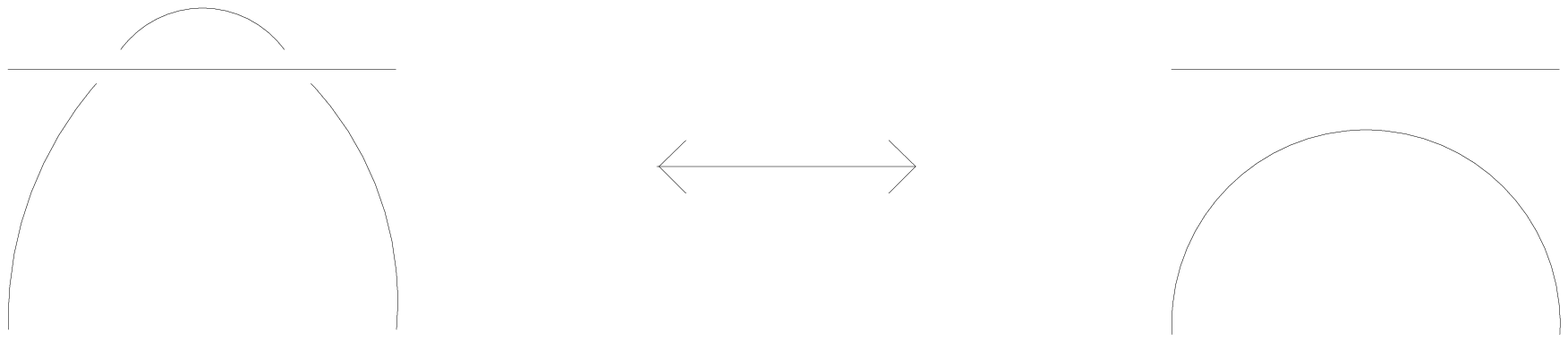} \end {figure}} \par  
{\begin {figure}[h] \hspace {1.6 cm} {\bf III)} \hspace {3 cm} 
\epsfysize = 0.9 cm 
\epsffile {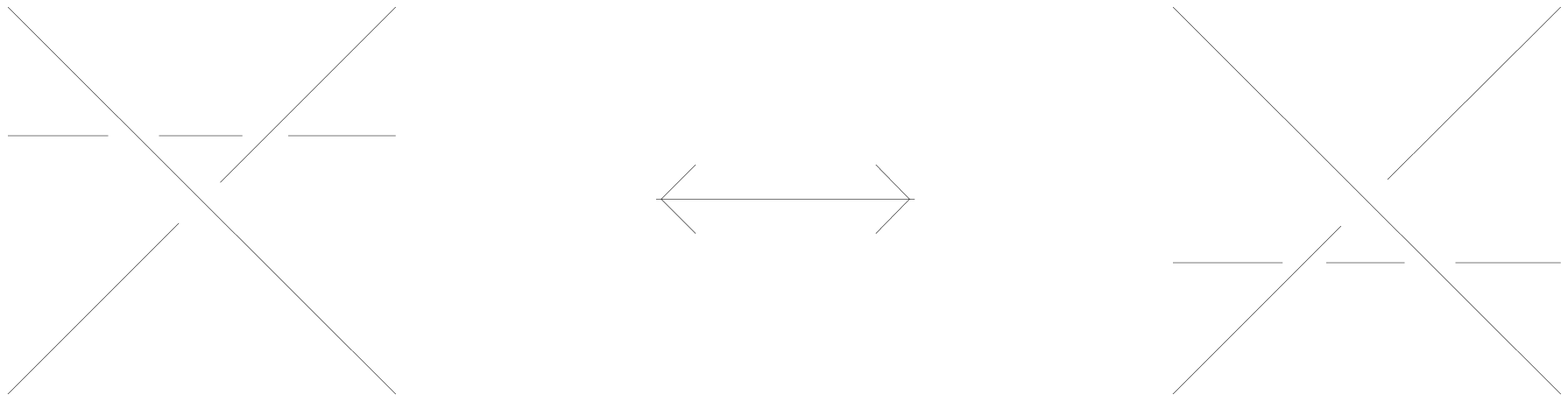} \end {figure}}  
\par 
These moves affect the set 
     that denotes the projection as follows. \par      
\medskip $\underline {First \quad move}$ : A pair $(i,i+1)$ or $(i+1,i)$ is 
     added or removed, while any number $j$ not less than $i$ is replaced by 
     $j+2$. \par \medskip  
     $\underline {Second \quad move}$: Two pairs $(i,j),(i+1,j+1)$ 
     or $(i+1,j),(i,j+1)$ are added or removed. Any number $k$ not less than  
${\rm min}(i,j)$ and less than ${\rm max}(i,j)$ is replaced by $k+2$, while any 
     number 
     $k$ not less than ${\rm max}(i,j)$ is replaced by $k+4$. \par \medskip
     $\underline {Third \quad move}$: Pairs $(i,j),(i',k),(j',k')$ where
     $|i'-i|=|j'-j|=|k'-k|=1$, 
     are replaced by $(i,k'),(i',j'),(j,k)$. \par \medskip      
     Sets connected through such moves denote ambient 
     isotopic knots and are thus identified. In addition, two sets $S_1$ and  
     $S_2$ are 
     identified if $\exists \epsilon \in \{ -1,1 \}$ and $c \in 
     \{ 1,2,...,2n \}$ such that $(i,j) \in S_1 \Rightarrow \Big( c+ \epsilon i
     \quad \rm {mod} (2n),c+ \epsilon  \quad \rm {mod} (2n) \Big) \in S_2$, 
     since these two sets denote the same 
     projection. The parameter $c$ indicates a change in the origin, while 
     $\epsilon = -1$
     indicates an orientation reversal. Finally, in order to identify mirror 
     symmetric knots, we consider equivalent two sets if $(i,j) \in S_1
     \Rightarrow (j,i) \in S_2$. \par \bigskip 
     \centerline {$\underline {\bf Knot \quad Invariants}$} \par \bigskip 
     Given a regular knot projection, one may inductively calculate any skein 
     invariant by repeatedly applying the formula \eject  
\par \smallskip {\begin {figure}[h]  
$$A f(\epsfysize=0.5 cm \epsffile {japfig1.eps})  
+ B f(\epsfysize=0.5 cm \epsffile {japfig2.eps})  
= C f(\epsfysize=0.5 cm \epsffile {japfig3.eps}) $$  
\end {figure}}   
\par \smallskip      
Eventually one ends up with a sum of 
multiples of the invariant of the trivial knot. This invariant is usually set 
to be equal to 1. When working on the set notation however, the following 
complication arises. While switching a crossing merely changes the type of the 
knot, eliminating a crossing converts a knot into a two-component link. During 
the inductive process one obtains a number of multicomponent links, which 
cannot be denoted through the Dowker notation as defined earlier. A way to 
overcome this problem is to use the unordered pair $\{ i,j \} $ each time a 
crossing 
$(i,j)$ or $(j,i)$ is removed. Therefore when an unordered pair appears in some 
set, it indicates the removal of a crossing from the original knot. It is thus 
possible to keep track of the emerging links, until one ends up with unlinks 
consisting of trivial knots. Then one continues by noticing that \par \smallskip
$$f(n \quad unlinked \quad circles) \quad = \quad \Big( {A+B \over c} \Big)
^{n-1}$$  
      In contrast to the previous notations, one need not stop here. 
      Inequivalent knots with identical skein invariants may be distinguished 
      through ``colorization" invariants. These invariants count the number of 
      distinct mappings from the projection {\it strands} (segments joining 
 undercrossings) to the set $\{ 1,2,...,k \}$, such that at each crossing the 
     three strands involved satisfy certain constraints (Ref. 12). The 
     constraints are defined in such a way so that if projection  becomes  
     through the application of a Reidemeister move, to each acceptable 
mapping of $P_1$ corresponds exactly one acceptable mapping of $P_2$. This 
     condition is necessary in order for the number of mappings to be an 
 invariant. In particular, let a strand $j$ separate the strand $i$ from the 
     strand $i+1$ as shown below, and let these strands be mapped to the 
``colors" $c_i$, $c_{i+1}$ and $c_j$. In order for the mapping to be allowed,
the strands must satisfy a relation $c_{i+1}=m_{c_ic_j}$, 
where $M_{ab}$ is a $k \times k$ matrix that defines the 
     colorization invariant. \par \medskip   
{\begin {figure}[h] \centerline {\epsfysize = 4 cm \epsffile {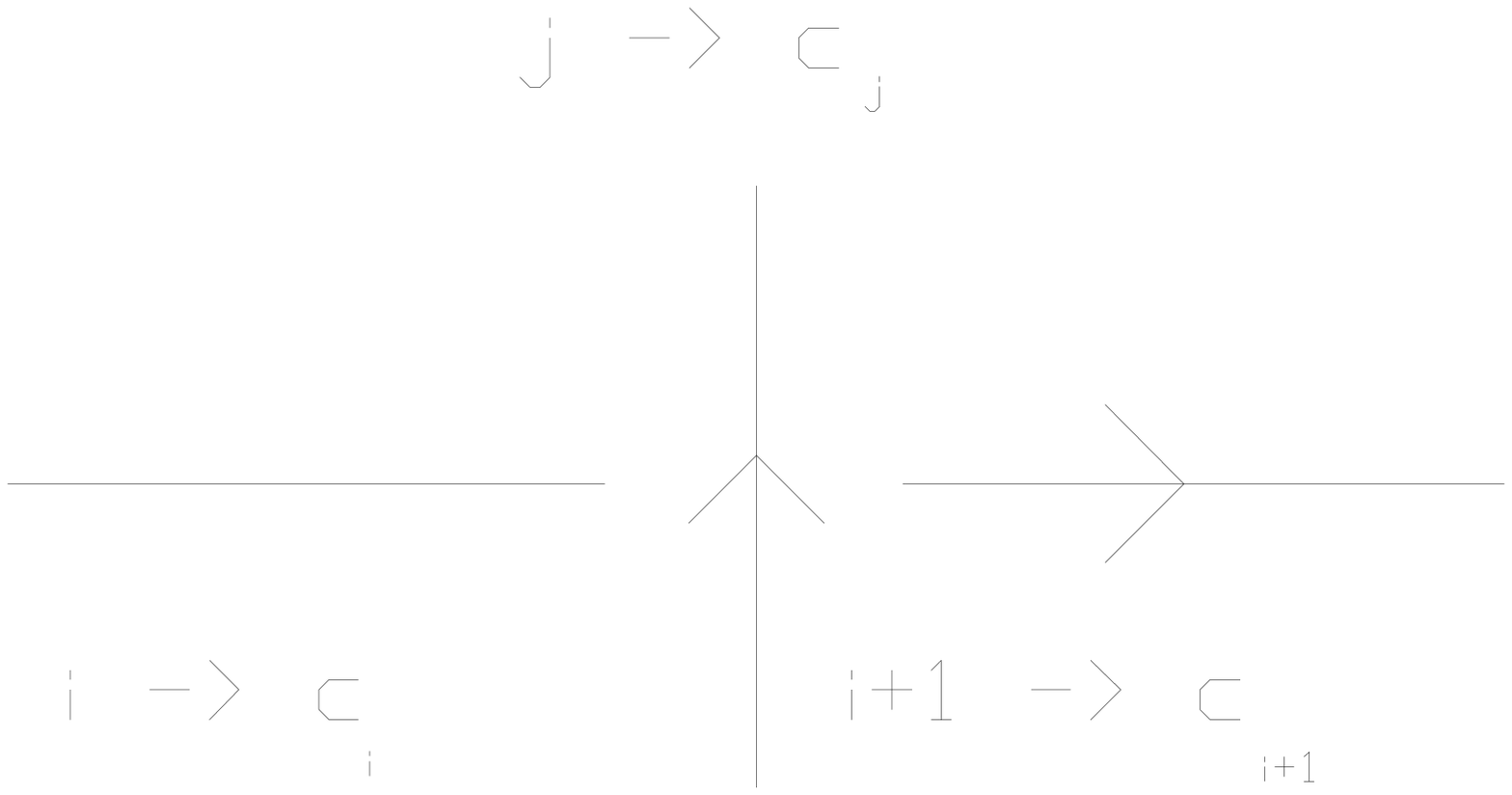}}  
    \end {figure}}  
     \par \medskip In order for the number of acceptable mappings 
     to be a true invariant, the matrix  has to fulfill the following 
     conditions. \par \medskip    
     1) $M_{rr}=r \quad \forall \quad r$ \par \smallskip           
     2) $M_{sr}=M_{tr} \quad \Leftrightarrow \quad s=t$ \par \smallskip      
     3) $M_{\alpha \beta} = \gamma, M_{\delta \alpha } = \epsilon, 
     M_{\gamma \beta } = \eta \quad \Rightarrow M_{\eta \gamma } =
     M_{\epsilon \beta }$ \par 
     \medskip Due to the enormous number of 
     possible such matrices, it is not feasible to check them one after the 
     other in order to obtain the necessary invariants. In practice we have 
     been able to obtain all matrices up to only 11 ``colors", which are 
     insufficient to distinguish knots with more than 7 crossings. Instead, 
     one uses special matrices that are guaranteed to obey the above 
     relations. One such category consists of matrices $M_{\alpha \beta } = 
     t \alpha + (1-t) \beta \quad {\rm mod} q$ where the integer $q$ is 
     prime with respect to both $t$ and $t-1$. (Corollary: $q$ is odd) 
     \par \medskip 	
     In order for the number of acceptable mappings to exceed the number of 
 ``colors", one obtains a homogeneous linear system with $n-1$ equations and 
$n-1$ unknowns, where $n$ is the crossing number. This system has non-zero 
     solutions if and only if its determinant is $0$; this determinant is in 
fact the {\bf Alexander} polynomial, and thus one may start by calculating this 
     polynomial before actually obtaining the colorization invariants.
     \par \medskip 
     A second category consists of matrices 
     $M_{g_\alpha g_\beta } = g_{\beta} \circ g_{\alpha } 
     \circ g_{\beta } ^{-1}$, 
where the ``colors" $g_{\rho }$ are elements 
     of a conjugacy class of some group. In particular, matrices derived from 
     permutation groups are extremely helpful.     
     \par \medskip One may easily check that 
     matrices of both these categories obey all necessary conditions to define 
     invariants. \par \bigskip 
     \centerline {$\underline {\bf Connected \quad Sums}$}  
     \bigskip When compiling tables of knots, one usually 
avoids {\it connected sums} and thus tabulates only {\bf prime} knots. There is
     a number of reasons for this. First, once prime knots have been tabulated, 
     one may obtain connected sums by combining prime knots. It is therefore 
     unnecessary to consider connected sums in the computer program before the 
     end, and one thus saves substantially in running time and computer 
     memory. Second, connected sums of the same prime knots possess identical 
     knot groups and cannot be distinguished through colorization invariants, 
     even if the connected sums are distinct. Third, if one uses the Dowker 
     notation, a set denoting the connected sum of $k$ prime knots may yield up 
     to $2^{k-1}$ distinct knots, and thus the notation becomes ambiguous.
\par \medskip 
     When one uses the braid word notation, or eventually arrives to the braid 
     word after using the closed lattice notation, one may check whether a 
     knot is prime or a connected sum as follows. If some braid generator 
     $\sigma _i \in B_n$, where $1 \leq i \leq n-1$, 
     appears exactly once in the braid word, the corresponding knot is 
     the connected sum of the closures of the two braids that are formed from 
     the products of the factors $\sigma _j$ for $j < i$ and $\sigma _k$ for 
     $k > i$. If both of these closures 
     are non-trivial knots, the initial knot is indeed a connected sum. The 
     problem however with this criterion is that if a generator appears 
     exactly once in the braid word, the knot is not presented in a minimal 
     form. Therefore one has to add at least one crossing and also increase 
     the braid index by one, in order to show that a knot is a connected sum.
     \par \medskip 
     This problem does not arise for the Dowker notation, which is one 
     additional reason that makes such a notation preferrable. If there is 
     some number $k$ between $1$ and $2n$ such that for any pair $(i,j)$ 
     belonging to 
     the set that represents the knot, 
     $(k-1) \cdot (k-j) \geq 0$, the knot is the connected sum of  two 
     sets, one consisting of pairs with numbers smaller than $k$, the other of 
     pairs with numbers larger than $k$. If both sets denote non-trivial knots, 
     the initial knot is a connected sum. One thus does not need to increase 
     the crossing number in order to show that a knot is a connected sum.
\eject 
\centerline {\bf 4. \quad COMPUTER PROGRAM RESULTS} \par \bigskip
     Using the Dowker notation we developed and ran a computer program in 
     order to tabulate knots. The main ideas of the algorithm were the 
     following. \par \medskip 
     We first obtained all regular projections whose crossing number does not 
     exceed a maximum value $n$. To do so, we started by considering all 
permutations in $S_k$ and all $k$-digit binary numbers, for $1 \leq k \leq n$. 
Let $S$ be the set 
     denoting a projection.  
     $p_i=j, b_i=0$ indicates that $(2i-1,2j) \in S$, while 
     $p_i=j, b_i=1$ indicates that $(2j,2i-1) \in S$. We then 
     eliminated the ``undrawable" sets and the sets denoting projections that 
     had already been recorded. \par    
     \medskip Once the regular projections were obtained, 
     we applied on each one of them all Reidemeister moves that do not 
     increase the crossing number. If a projection was not connected to one 
     checked before, the projection was stored in the computer's memory. If it 
     was connected to two or more projections already stored in the memory, 
     these projections would be identified. When the process was completed, 
     the computer recorded all projections stored in the memory that were not 
     identified later to any previous projections. \par 
\medskip 
When $n$ was set $14$, the number of recorded projections as a function of the 
    crossing number came out as follows. (Connected sums are not included; 
    mirror symmetric and orientation reversed knots are listed exactly once, 
    even if they are not ambient isotopic) \par \medskip 
\centerline {\vbox {\def\*{\hphantom{0}} \offinterlineskip 
\halign {\strut#&\vrule#\quad&\hfil#\hfil&\quad\vrule#\quad&\hfil#\hfil&\quad
\vrule#\cr \noalign {\hrule}
&&\omit \it Number of Crossings && \it Recorded projections &\cr \noalign 
{\hrule } 
&&\*0&& \*\*\*\*1&\cr 
&&\*1&& \*\*\*\*0&\cr 
&&\*2&& \*\*\*\*0&\cr 
&&\*3&& \*\*\*\*1&\cr 
&&\*4&& \*\*\*\*1&\cr 
&&\*5&& \*\*\*\*2&\cr 
&&\*6&& \*\*\*\*3&\cr 
&&\*7&& \*\*\*\*7&\cr 
&&\*8&& \*\*\*21&\cr 
&&\*9&& \*\*\*49&\cr 
&&10&& \*\*165&\cr 
&&11&& \*\*552&\cr 
&&12&& \*2191&\cr 
&&13&& 29781&\cr \noalign 
{\hrule}}}} \par \medskip
    All other projections are equivalent to one or more of the recorded 
    projections. It is also possible that some of the projections of the table 
    above are ambient isotopic to each other, but may only be connected 
    through Reidemeister moves involving $15$ or more crossings. In order to 
    obtain the distinct knot types we calculated first the Alexander 
    polynomials of these projections. For projections having identical 
    Alexander polynomials we calculated colorization invariants for 
permutation groups $S_p$, where $1 \leq p \leq m$ and $m$ is a fixed parameter. 
For $m=7$ we were 
able to show that all knots up to $11$ crossings are inequivalent, while (at 
least) $2176$ out of the $2191$ knots with crossing number $12$ are 
inequivalent 
    to each other and to knots with smaller crossing numbers (for detailed 
    tables and explicit proofs, see Ref. 13).  These results are consistent to 
    results recently reported by J. Hoste, M. Thistlethwaite and J. Weeks in 
    Ref. 1. \par 
\medskip 
     One observes that the whole program needs merely two input numbers $n$ and 
     $m$, which indicate the maximum crossing number of the projections 
     considered, and the largest permutation group used to obtain colorisation 
     invariants. As stated earlier, we eventually set these parameters equal 
     to $n=14$ and $m=7$. \par 
\medskip 
     While running the program for smaller parameters, we also noticed the 
     following. First, in order to establish all ambient isotopies for knots 
     up to $c$ crossings, one had to set  
     $n \geq {3c \over 2} - 3$ for $5 \leq c \leq 11$. To be exact, the table 
     obtained 
     was the following. \par \medskip 
\centerline {\vbox {\def\*{\hphantom{0}} \offinterlineskip 
\halign {\strut#&\vrule#\quad&\hfil#\hfil&\quad\vrule#\quad&\hfil#\hfil&\quad
\vrule#\cr \noalign {\hrule}
&&\omit \it c && \it n &\cr \noalign 
{\hrule } 
&&\*0&& \*0&\cr 
&&\*1&& \*0&\cr 
&&\*2&& \*0&\cr 
&&\*3&& \*3&\cr 
&&\*4&& \*4&\cr 
&&\*5&& \*5&\cr 
&&\*6&& \*6&\cr 
&&\*7&& \*8&\cr 
&&\*8&& \*9&\cr 
&&\*9&& 11&\cr 
&&10&& 12&\cr 
&&11&& 14&\cr \noalign  
{\hrule}}}} \par \medskip
     Second, in order to prove that the tabulated knots were distinct, for
     $c \leq 8$
 the Alexander polynomials were sufficient. For $c=9$ or $10$ one had to set 
 $m \geq 5$, 
     while for $c=11$ one had to set $m \geq 7$. \par \bigskip 
     \centerline {ACKNOWLEDGMENTS} \par \bigskip      
     It is a pleasure to 
     acknowledge Wei Chen, K. Anagnostopoulos, J. Przytycki, S. Lambropoulou, 
     S. Garoufalidis, C. Shubert and D. Westbury for their help and/or 
     interesting discussions. I would also like to thank the Brookhaven 
     National Laboratory, the Weizmann Institute of Science and the University 
     of Gottingen for their hospitality during parts of the work. Finally I 
     would like to thank the Waseda University of Tokyo for organizing the 
 conference {\it Knots 96} and for giving me the opportunity to present this 
     project during the conference. \par \bigskip \centerline 
{\bf REFERENCES} \par \bigskip 
1) M.B. Thistlethwaite, {\it L.M.S. Lecture Notes} {\bf 93}, 1-76, Cambridge 
University Press, 1985; J. Hoste, M. Thistlethwaite, J. Weeks, {\it Knot 
	Tabulation Progress Report}, July 1996.
\par \medskip 
2) J.W. Alexander, G.B. Briggs, {\it Annals of Mathematics} {\bf 2}, vol. 28, 
563-586, (1927); H. Brunn, {\it Verh. des intern. Math. Congr.} {\bf 1}, 
256-259, (1897); H.R. Morton, {\it Math. Proc. Cambridge Philos. Soc.} 
{\bf 99}, 247-260 
(1986); P. Vogel, {\it Comment. Math. Helvetici} {\bf 65}, 104-113 (1990); 
S. Yamada, 
     {\it Invent. Math.} {\bf 89}, 347-356 (1987). 
\par \medskip 
3) S. Lambropoulou, C. Rourke, {\it Topology and its Applications} (to be 
	published).
\par \medskip 
4) A.A. Markov, {\it Recusil Mathematique Moscou} {\bf 1}, (1935); J.S. Birman, 
{\it Ann. of Math. Stud.} {\bf 82}, Princeton University Press, Princeton, 
1974; 
D. Bennequin, {\it Asterisque} {\bf 107-108}, 87-161 (1993); P. Traczyk, 1992 
	(preprint).
\par \medskip
5) J.W. Alexander, {\it Trans. Amer. Math. Soc.} {\bf 30}, 275-306 (1926); 
V.F.R. Jones, {\it Notices of AMS} {\bf 33}, 219-225 (1986); P. Freyd, 
D. Yetter, 
      J. Hoste, W.B.R. Lickorish, K.C. Millet, A. Ocneanu {\it Bull. Ams} 
{\bf 12}, 239-246 (1985); J.H. Przytycki, K.P. Traczyk, Kobe J. Math. 4, 115 
      (1987).
\par \medskip 
       6) D. Rolfsen, {\it Knots and Links}, Berkeley, 1976.
\par \medskip 
7) N. Pippenger, {\it Discrete Appl. Math.} {\bf 25}, 273-278 (1989); 
E.J. Janse 
van Rensburg, {\it J. Phys.} {\bf A 25}, 1031-1042 (1992); 
E.J. Janse van Rensburg, 
S.G. Whittington, {\it J. Phys.} {\bf A 24}, 5553-5567; 
C.E. Soteros, D.W. Sumners, 
S.G. Whittington, {\it Math. Proc. Camb. Phil. Soc.} {\bf 111}, 75-91 (1992); 
D.W. Sumners, S.G. Whittington, {\it J. Phys.} {\bf A 21}, 1689-1694 (1988).
\par \medskip 
8) C. Aneziris, {\it The Braid Word of a closed Self-Avoiding Walk}, (under 
       preparation).
\par \medskip 
 9) Y.A. Diao, {\it Journal of Knot Theory and its Ramifications} {\bf 2}, 
413-427 (1993).  
\par \medskip 
10) C. Jourdan, {\it Course d' Analyse}, 1893; O. Veblen, {\it Trans. Amer. 
Math. 
     Soc.} {\bf 6}, 83 (1905).
\par \medskip 
11) K. Reidemeister {\it Abh. Math. Sem. Univ. Hamburg} {\bf 5}, 7-23 (1927).
\par \medskip 
12) R.H. Fox, {\it Canadian J. Math.} {\bf XXII(2)}, 193-201 (1970); R. Fenn, 
C. Rourke, {\it Journal of Knot Theory and its Ramifications} {\bf 1}, 
343-406 (1992).
\par \medskip 
     13) C. Aneziris, http://sgi.ifh.de/$\sim$aneziris/contents.html/   
\end{document}